\documentclass[11pt]{amsart}
\usepackage{amsmath,amssymb}
\usepackage[hyperindex]{hyperref} 
\usepackage{amsrefs,booktabs,verbatim}
\usepackage{xcolor}

\title{Questions about determinants and polynomials}
\author{Steve Fisk}
\date{\today}
 
\newtheorem{question}{Question}
\theoremstyle{definition}
\newtheorem{definition}{Definition}
\newtheorem{remark}{Remark}
\newtheorem{example}{Example}

\newenvironment{aside}{}{}

\newcommand{\ulace}{\ensuremath{\stackrel{U}{\longleftarrow}}}

\newcommand{\stabled}[1]{{\mathcal{H}_{#1}}} 

\newcommand{\polypos}[1]{P^{+}_{#1}}

\newcommand{\rup}[1]{{\text{\texttt{U}}_{#1}}}

\begin{document}
\maketitle

This article is an exposition of some  questions about
determinants and polynomials. We are interested in the following
classes of polynomials:

\begin{definition}\ 


  $\rup{d} = \left\{\text{\parbox{3.5in}{All polynomials
        $f(x_1,\dots,x_d)$ with real coefficients such that
        $f(\sigma_1,\dots,\sigma_d)\ne0$ for all
        $\sigma_1,\dots,\sigma_d$ in the open upper half plane. 
        }}\right.$

 $\polypos{d} = $ all polynomials in $\rup{d}$ with all positive coefficients.

$\stabled{1} = $ all polynomials with all roots in the closed left half
plane. 
\end{definition}

Polynomials in $\stabled{1}$ are called \emph{stable polynomials} or
sometimes \emph{Hurwitz stable} polynomials. The polynomials in
$\rup{d}$ are called \emph{upper polynomials}.  $\polypos{1}$ consists
of polynomials with all negative roots, and so are also stable.
Surveys of these classes of polynomials can be found in
\cite{fisk-upper} and \cite{fisk-stable}. We say $f(x)$ and $g(x)$
\emph{interlace}, written $f \ulace g$, if $f+y g\in\rup{2}$. This is
equivalent to the usual definition that the roots of $f$ and $g$
alternate.

\begin{question}\label{ques-1}
  Suppose that $f = \sum_0^n a_iy^i$ is in $\polypos{1}$. Form the polynomial

\begin{equation}\label{eqn:ques-1}
F(x)=
\begin{vmatrix} a_0 & a_1 \\ 0 & a_0 \end{vmatrix} +
\begin{vmatrix} a_1 & a_2 \\ a_0 & a_1 \end{vmatrix}\,x +
\cdots +
\begin{vmatrix} a_n & 0 \\ a_{n-1} & a_n \end{vmatrix}\, x^n
\end{equation}

Show that this is in $\polypos{1}$.

\end{question}

\begin{remark}
  See \cite{sagan} for some recent results on this problem.
  
  Using the criteria for stability for polynomials of small degree \cite{fisk} we
  can show that the polynomials of degree at
  most 4 are stable. In the case of degree two it follows easily that the
  polynomials are actually in $\polypos{1}$.

\begin{description}
\item[degree $=2$]  If we write $f(x)=(x+a)(x+b)$ where $a,b$
  are positive then

\[
F(x)= 
a^2 b^2 + \bigl(a^2+b a+b^2\bigr)x + x^2
\]
Since all coefficients are positive $F(x)$ is stable. The discriminant
is $\left(a^2-b a+b^2\right) \left(a^2+3 b a+b^2\right)$ which is
positive for positive $a,b$, so $F$ is in $\polypos{1}$.
\item[degree $=3$]  Write $f(x)=(x+a)(x+b)(x+c)$ with positive
  coefficients. We find that
  \begin{align*}
    F(x) &= a^2 b^2 c^2 +
\bigl( a^2 b^2 + a^2 b c + a b^2 c + a^2 c^2 + a b c^2 + b^2
c^2\bigr)x \\& +
 \bigl(a^2 + a b + b^2 + a c + b c + c^2\bigr)x^2 + x^3
  \end{align*}
  If we write $F = \alpha_0 + \alpha_1 x + \alpha_2 x^2 + x^3$ then
  the criterion to be stable is $\alpha_1\alpha_2-\alpha_0>0$. If we
  compute this expression we get a sum of 19 monomials, all with
  positive coefficients, so $F$ is stable.
\item[degree $=4$] 
  In this case we compute the criterion  to be stable, and it is a sum
  of 201 monomials, all with positive coefficients.

\end{description}

\end{remark}

\begin{question}\label{ques-2}
     Generalize Question~\ref{ques-1} by considering

\[T_k(f) = \sum_i \begin{vmatrix}a_i & a_{i+k} \\ a_{i-k} &
  a_i \end{vmatrix}x^i\] 
Show that $T_k(f)$ is in $\polypos{1}$ if $f\in\polypos{1}$.
\end{question}

\begin{remark}
A computer algebra calculation shows
that $T_2(f)$ is stable when $f$ has degree four. If $f$
has degree $n$ and $k> n/2$ then
\[ T_k(f) = \sum a_i^2 x^i = f\ast f\]
where $\ast$ is the Hadamard product. Thus, $T_k(f)\in\polypos{1}$ if
$k>n/2$.  

Computation shows that $T_k(f)$ and $T_j(f)$ do not generally
interlace. However, it does appear that there is a $g\in\polypos{1}$
that interlaces every $T_k(f)$. 
\end{remark}

\begin{question}\label{ques-2a}
Choose a positive
   integer $d$, and let $\sum_0^n a_i\,x^i\in\polypos{1}$. Form
   \begin{equation}
     \label{eqn-2a}
     F(x) = 
\sum_i
x^i\, 
     \begin{vmatrix}
       a_i & \dots & a_{i+d} \\ a_{i-1} & \dots & a_{i+d-1} \\ \vdots & & \vdots 
       \\ a_{i-d} & \dots & a_i
     \end{vmatrix}
   \end{equation}
   
Show that $F(x)$ is in $\polypos{1}$.
\end{question}

\begin{remark}
  If $d=1$ this is just Question~\ref{ques-1}.
\end{remark}

\begin{question}\label{ques-2d}
  If $\sum_0^n a_ix^i\in\polypos{1}$ then construct the matrix
\[
M=\begin{pmatrix}
  a_0 & a_1 & a_2 & a_3 & \dots \\
  0 & a_0 & a_1 & a_2 & \dots \\
  0 & 0 & a_0 & a_1 & \dots \\
  \vdots & \vdots
\end{pmatrix}
\]
For any positive integer $d$ construct a new matrix $M'$ by replacing
each entry  of $M$ by the determinant of the $d$ by $d$ matrix whose
upper left corner is that entry. Show that $M'$ is totally positive.
\end{question}

\begin{remark}
  By the Aissen-Schoenberg theorem, if $M'$ is totally positive then
  the polynomial corresponding to the first row is in
  $\polypos{1}$. Thus this question implies Questions \ref{ques-1} and \ref{ques-2a}.
\end{remark}

\begin{question}\label{ques-2b}
  Suppose that $f\ulace g$ in $\polypos{1}$. If $f=\sum a_ix^i$, $g =
  \sum b_ix^i$ then show that
\[
\sum_i
\begin{vmatrix}
  a_i & a_{i+1} \\ b_i & b_{i+1}
\end{vmatrix}\,x^i \in \polypos{1}
\]
\end{question}
\begin{remark}
  This implies Question~\ref{ques-1}. Since $(x+1)f(x) \ulace f(x)$ we
  have 
\[
\sum_i
\begin{vmatrix}
  a_i + a_{i-1} & a_{i+1}+a_i \\ a_i & a_{i+1} 
\end{vmatrix}\,x^i =
\sum_i
\begin{vmatrix}
a_{i-1} & a_i \\ a_i & a_{i+1}  
\end{vmatrix} \,x^i
\]
which is Question~\ref{ques-1}.
\end{remark}

\begin{question}\label{ques-2c}
  Suppose $f = \sum a_{i,j}x^iy^j\in\polypos{2}$. For any positive
  integer $d$ show that
\[
\sum_i x^i\, 
\begin{vmatrix}
  a_{i,0} & \dots & a_{i,d} \\
\vdots && \vdots \\
a_{i+d,0} & \dots & a_{i+d,d}
\end{vmatrix}
\in\polypos{1}
\]
\end{question}
\begin{remark}
  Since $f\ulace g$ is equivalent to $f + y g\in\rup{2}$ we see
  that Question~\ref{ques-2c} implies Question~\ref{ques-2b}.
\end{remark}

\begin{question}\label{ques-3}
   Consider a polynomial $f$ in $\polypos{3}$ 

\[
\begin{matrix}
  f_{0,0}(x) &+  & f_{0,1}(x)y &+& f_{0,2}(x)y^2 &+& \dots \\[.2cm]
  f_{1,0}(x)z &+  & f_{1,1}(x)yz &+&  f_{1,2}(x)y^2z &+&\dots \\[.2cm]
  f_{2,0}(x)z^2 &+  & f_{2,1}(x)yz^2 &+ &  f_{2,2}(x)y^2z^2 &+& \dots \\[.2cm]
\end{matrix}
  \]

We construct a new polynomial by replacing each term by  the $k$ by
$k$ determinant based at that term. If $k=2$ then the polynomial is 
$F(x,y,z)=$
\[
\begin{matrix}
\begin{vmatrix} f_{0,0} & f_{0,1}\\ f_{1,0} & f_{1,1} \end{vmatrix}
& +&
\begin{vmatrix} f_{0,1} & f_{0,2}\\ f_{1,1} &  f_{1,2} \end{vmatrix}\, y
&+&
\begin{vmatrix} f_{0,2} & f_{0,2}\\ f_{1,1} & f_{1,2} \end{vmatrix}\, y^2
& +& \dots \\[.5cm]
\begin{vmatrix} f_{1,0} & f_{1,1}\\ f_{2,0} & f_{2,1} \end{vmatrix}\, z
& +& 
\begin{vmatrix} f_{1,1} & f_{1,2}\\ f_{2,1} & f_{2,2} \end{vmatrix}\, yz
& +& 
\begin{vmatrix} f_{1,2} & f_{1,3}\\ f_{2,2} & f_{2,3} \end{vmatrix}\,y^2 z
&+& \dots \\[.5cm] 
\end{matrix}
\]

Show that  if $k=2$ then $F(x,\alpha,\beta)$ is stable for positive
$\alpha,\beta$. In addition, for all $k$ all coefficients have the same sign.

\end{question}

\begin{remark}

  If all $f_{i,j}(x)$ are constant, so that $F(x,y,z)=G(y,z)$ where
  $G\in\polypos{2}$, then $G(x,0)$ is equation \eqref{ques-2a}, and so would be
  in $\polypos{1}$, rather than just a stable polynomial.

  Here are two simple cases. If we take $f = (x+y+z)^2$ then 
\[ F(x,y,z) = -2(x^2+y+z) \]
This is not stable, but is stable for positive $y,z$.

If $f= (x+y+z)^3$ then
\[F(x,y,z) = -3 \left(x^4+3 y x^2+3 z x^2+y^2+z^2+3 y z\right) \]
Again, $F(x,y,z)$ is not stable, and can be checked to be stable for
positive $y,z$.

If we consider $k=3$ and $f=(x+y+z)^3$ then $F(x,y,z) = -9(x^3+y+z)$
which is not a stable polynomial and
$F(x,\alpha,\beta)\not\in\polypos{1}$ for all positive $\alpha,\beta$,
but it does have all negative coefficients.
\end{remark}

\begin{question}\label{ques-4a}
  Suppose that $f\in\polypos{1}$. Show that the determinant of the matrix below
  is stable. We know that it is positive for positive $x$.

\begin{equation}\label{eqn-4a}
F(x)=\begin{pmatrix}
  f & f' & \dots & f^{(d)}\\
\vdots &\vdots&& \vdots \\
  f^{(d)} & f^{(d+1)} & \dots & f^{(2d)}
\end{pmatrix}
\end{equation}

\end{question}

\begin{question}\label{ques-4}
     Suppose $f\in\polypos{2}$, and write $f = \sum f_i(x)y^i$.
     Consider the polynomial 
\begin{equation}\label{eqn-4}
F(x)=\begin{vmatrix}
  f_i & f_{i+1} & \dots & f_{i+d}\\
  f_{i+1} & f_{i+2} & \dots & f_{i+d+1} \\
\vdots &\vdots&& \vdots \\
  f_{i+d} & f_{i+d+1} & \dots & f_{2d}\\
\end{vmatrix}
\end{equation}

Show that $F(x)$ is stable for all positive integers $d$ and
non-negative integers $i$.

\end{question}

\begin{remark}
  We know that this holds for $d=1$. If we consider $i=0$ and  $f(x+y)$ where
  $f\in\polypos{1}$ then from \eqref{eqn-4} we have an assertion that appears to be
  stronger than Question~\ref{ques-4a}:

\begin{equation*}
\begin{vmatrix}
  f & f' &f^{(2)}/2!& \dots & f^{(d)}/d!\\
\vdots &\vdots&\vdots&\vdots&  \vdots \\
  f^{(d)}/d! & f^{(d+1)}/(d+1)! &\dots &  \dots & f^{(2d)}/(2d)!
\end{vmatrix} \qquad \text{is stable.}
\end{equation*}

\end{remark}

\begin{question}\label{ques-5}
     Suppose $f\in\polypos{2}$, and write $f = \sum f_i(x)y^i$. Form
     the infinite matrix
\begin{equation}\label{eqn-5}
\begin{pmatrix}
  f_0 & f_1 & f_2 & f_3 & \dots \\
  0 & f_0 & f_1 & f_2  & \dots \\
  0 & 0 & f_0 & f_1  & \dots \\
  0 & 0 & 0  & f_0  & \dots \\
  \vdots &\vdots&\vdots&\vdots&\ddots
\end{pmatrix}
\end{equation}
Show that the resulting matrix is totally stable. That is, every
minor is a stable polynomial.
\end{question}
\begin{remark}
  If we substitute a positive value for $x$ then the resulting matrix
  is well known to be totally positive. Question~\ref{ques-4} is a
  consequence of this question since the the determinant in
  \eqref{eqn-4} is a minor of \eqref{eqn-5}.

  For example, if $f = (x+y)^n$ then it appears that all minors of
  \eqref{eqn-5} have the form $cx^s$ for positive $c$ and non-negative
  integer $s$.

\end{remark}

\begin{question} \label{ques-7} 
  Suppose that $f = \prod(x+ b_i y + c_i)$ where $b_i,c_i$ are
  positive, and   write $f = \sum a_{ij}x^iy^j$. Consider the matrix

\[
\begin{pmatrix}
  a_{0d} & \dots & a_{00} \\
\vdots & & \vdots \\
a_{dd} & \dots & a_{d0}
\end{pmatrix}
\]

Show that the matrix is totally positive for any positive integer $d$.
 \end{question}

 \begin{remark}
   $f$ is in $\polypos{2}$, but the assertion fails for arbitrary
   polynomials in $\polypos{2}$. However, since consecutive rows are
   the coefficients of interlacing polynomials, all two by two
   determinants are positive for any $f\in\polypos{2}$.

   If we take $\prod_1^3(x+b_iy+c_i)$ and  $d=2$ then the determinant
   is a sum of 7 monomials with all positive coefficients. 

   Here's an example where it fails for $f\in\polypos{2}$.     \begin{gather*}
      M =
      \begin{pmatrix}
        1&0&0\\0&1&0\\0&0&1
      \end{pmatrix}
+
x\,\begin{pmatrix}
  13& 9& 7  \\9 & 7& 5 \\7& 5& 4
\end{pmatrix}
+
y\,
\begin{pmatrix}
  5& 7& 8 \\7& 11& 12\\8& 12& 14
\end{pmatrix}\\
\intertext{The three matrices are positive definite, so the
  determinant of $M$ is in $\polypos{2}$, and equals} 1 + 24\,x +
16\,x^2 + 2\,x^3 + 30\,y + 164\,x\,y +
62\,x^2\,y + 22\,y^2 + 64\,x\,y^2 + 4\,y^3 \\
\intertext{with coefficient matrix}
\begin{pmatrix}
 4&0&0&0\\22&64&0&0\\30&164&62&0\\ 1&24&16&2 
\end{pmatrix}
    \end{gather*}

    The determinant of the three by three matrix in the lower left
    corner is $-1760$. Of course, all the two by two submatrices have
    positive determinant.

 \end{remark}

 \begin{question}
   When is the product of two totally stable matrices a totally stable matrix?
 \end{question}

 \begin{question}
   A \emph{totally upper matrix} has the property that every minor is
   either zero or a polynomial with all real roots.
   \begin{enumerate}
   \item What are constructions of totally upper matrices and totally
     stable matrices?
   \item When is the product of two totally upper matrices a totally
     upper matrix?
   \end{enumerate}
 \end{question}

\end{document}